\documentclass[12pt]{article}
\usepackage{latexsym,amssymb,amsmath,enumerate,geometry,float,cite,tikz,setspace}
\newtheorem{definition}{Definition}
\newtheorem{theorem}[definition]{Theorem}
\newtheorem{corollary}[definition]{Corollary}

\usepackage{lineno}

\begin{document}


\onehalfspace

\title{Relating $2$-rainbow domination\\ to weak Roman domination}

\author{Jos\'{e} D. Alvarado$^1$, Simone Dantas$^1$, and Dieter Rautenbach$^2$}

\date{}

\maketitle

\begin{center}
{\small 
$^1$ Instituto de Matem\'{a}tica e Estat\'{i}stica, Universidade Federal Fluminense, Niter\'{o}i, Brazil\\
\texttt{josealvarado.mat17@gmail.com, sdantas@im.uff.br}\\[3mm]
$^2$ Institute of Optimization and Operations Research, Ulm University, Ulm, Germany\\
\texttt{dieter.rautenbach@uni-ulm.de}
}
\end{center}

\begin{abstract}
Addressing a problem posed by Chellali, Haynes, and Hedetniemi (Discrete Appl. Math. 178 (2014) 27-32)
we prove $\gamma_{r2}(G)\leq 2\gamma_r(G)$ for every graph $G$,
where $\gamma_{r2}(G)$ and $\gamma_r(G)$ denote the 
$2$-rainbow domination number and the weak Roman domination number of $G$, 
respectively.
We characterize the extremal graphs for this inequality that are $\{ K_4,K_4-e\}$-free, 
and show that the recognition of the $K_5$-free extremal graphs is NP-hard.
\end{abstract}

{\small 

\medskip


\medskip


}

\section{Introduction}\label{section1}

We consider finite, simple, and undirected graphs, and use standard terminology and notation. 

Rainbow domination was introduced in \cite{bhr}.
Here we consider the special case of $2$-rainbow domination. 
A {\it $2$-rainbow dominating function} of a graph $G$ is a function $f:V(G)\to 2^{\{ 1,2\}}$ 
such that $\bigcup_{v\in N_G(u)}f(v)=\{ 1,2\}$ for every vertex $u$ of $G$ with $f(u)=\emptyset$.
The {\it weight of $f$} is $\sum_{u\in V(G)}|f(u)|$.
The {\it $2$-rainbow domination number $\gamma_{r2}(G)$ of $G$} 
is the minimum weight of a $2$-rainbow dominating function of $G$,
and a $2$-rainbow dominating function of weight $\gamma_{r2}(G)$ is {\it minimum}.
Weak Roman domination was introduced in \cite{hh}.
For a graph $G$, a function $g:V(G)\to \mathbb{R}$, 
and two distinct vertices $u$ and $v$ of $G$, let
$$g_{v\to u}:V(G)\to \mathbb{R}:
x\mapsto
\left\{
\begin{array}{rl}
g(u)+1 &, x=u,\\
g(v)-1 &, x=v,\mbox{ and}\\
g(x) &, x\in V(G)\setminus \{ u,v\}.
\end{array}
\right.
$$
A set $D$ of vertices of $G$ is {\it dominating} if every vertex in $V(G)\setminus D$ has a neighbor in $D$.
A {\it weak Roman dominating function} of $G$ is a function $g:V(G)\to \{ 0,1,2\}$ 
such that every vertex $u$ of $G$ with $g(u)=0$ 
has a neighbor $v$ with $g(v)\geq 1$
such that the set $\{ x\in V(G):g_{v\to u}(x)\geq 1\}$ is dominating.
The {\it weight of $g$} is $\sum_{u\in V(G)}g(u)$.
The {\it weak Roman domination number $\gamma_r(G)$} of $G$
is the minimum weight of a weak Roman dominating function of $G$,
and a weak Roman dominating function of weight $\gamma_r(G)$ is {\it minimum}.

For a positive integer $k$, let $[k]$ be the set of positive integers at most $k$.

\bigskip

\noindent In \cite{chh} Chellali, Haynes, and Hedetniemi show that 
$\gamma_r(G)\leq \gamma_{r2}(G)$ for every graph $G$,
and pose the problem to upper bound the ratio $\frac{\gamma_{r2}(G)}{\gamma_r(G)}$ (cf. Problem 17 in \cite{chh}).
In the present paper we address this problem.
As we shall see in Theorem \ref{theorem1} below, 
$\frac{\gamma_{r2}(G)}{\gamma_r(G)}\leq 2$ for every graph $G$.
While the proof of this inequality is very simple, the extremal graphs are surprisingly complex.
We collect some structural properties of these graphs in Theorem \ref{theorem1},
and characterize all $\{ K_4,K_4-e\}$-free extremal graphs in Corollary \ref{corollary2},
where $K_n$ denotes the complete graph of order $n$,
and $K_n-e$ arises by removing one edge from $K_n$.
In contrast to this characterization, 
we show in Theorem \ref{theorem2} 
that the recognition of the $K_5$-free extremal graphs is algorithmically hard,
which means that these graphs do not have a transparent structure.
In our last result, Theorem \ref{theorem3},
we consider graphs whose induced subgraphs are extremal.

The weak Roman domination number was introduced as a variant 
of the Roman domination number $\gamma_R(G)$ of a graph $G$ \cite{st}.
For results concerning the ratio $\frac{\gamma_{r2}(G)}{\gamma_R(G)}$ see \cite{cr,ff,wx}.

\section{Results}

\begin{theorem}\label{theorem1}
If $G$ is a graph, then $\gamma_{r2}(G)\leq 2\gamma_r(G)$.
Furthermore, if $\gamma_{r2}(G)=2\gamma_r(G)$ 
and $g:V(G)\to \{ 0,1,2\}$ is a minimum weak Roman dominating function of $G$, then
\begin{itemize}
\item there is no vertex $x$ of $G$ with $g(x)=2$, and
\item if $V_1=\{ v_1,\ldots,v_k\}$ is the set of vertices $x$ of $G$ with $g(x)=1$, 
then $V(G)\setminus V_1$ has a partition into $2k$ sets $P_1,\ldots,P_k,Q_1,\ldots,Q_k$ 
such that for every $i\in [k]$,
\begin{itemize}
\item $P_i=\{ u\in V(G)\setminus V_1:N_G(u)\cap V_1=\{ v_i\}\}$ is non-empty and complete for $i\in [k]$, and
\item every vertex in the possibly empty set $Q_i$ is adjacent to every vertex in $\{ v_i\}\cup P_i$.
\end{itemize}
\end{itemize}
\end{theorem}
{\it Proof:} Let $g:V(G)\to \{ 0,1,2\}$ is a minimum weak Roman dominating function of $G$.
Clearly, 
$f:V(G)\to 2^{[2]}$ with
$$f(x)=
\left\{
\begin{array}{rl}
\emptyset, & g(x)=0\mbox{ and }\\
\{ 1,2\}, & g(x)>0
\end{array}
\right.
$$
is a $2$-rainbow dominating function of $G$, 
which immediately implies
\begin{eqnarray}
\gamma_{r2}(G) & \leq & \sum_{u\in V(G)}|f(u)|\leq 2\sum_{u\in V(G)}g(u)=2\gamma_r(G).\label{e1}
\end{eqnarray}
Now, let $\gamma_{r2}(G)=2\gamma_r(G)$, which implies that equality holds throughout (\ref{e1}).
This implies that there is no vertex $x$ of $G$ with $g(x)=2$.
Let $V_1=\{ v_1,\ldots,v_k\}$ be the set of vertices $x$ of $G$ with $g(x)=1$.
For $i\in [k]$, let $P_i=\{ u\in V(G)\setminus V_1:N_G(u)\cap V_1=\{ v_i\}\}$,
that is, for $u\in P_i$, the only neighbor $v$ of $u$ with $g(v)\geq 1$ is $v_i$.
Therefore, the set $\{ x\in V(G):g_{v_i\to u}(x)\geq 1\}$ is dominating,
which implies that $P_i$ is complete.
If $P_i=\emptyset$ for some $i\in [k]$, then 
$f':V(G)\to 2^{[2]}$ with
$$f'(x)=
\left\{
\begin{array}{rl}
\emptyset, & g(x)=0,\\
\{ 1,2\}, & x\in V_1\setminus \{ v_i\},\mbox{ and }\\
\{ 1\}, & x=v_i
\end{array}
\right.
$$
is a $2$-rainbow dominating function of $G$ of weight less than $2\gamma_r(G)$,
which is a contradiction.
Hence, for every $i\in [k]$, the set $P_i$ is non-empty and complete.

For $u\in V(G)\setminus (V_1\cup P_1\cup\cdots\cup P_k)$,
let $i(u)$ be the smallest integer in $[k]$ such that $v_{i(u)}$ is a neighbor of $u$
and the set $\{ x\in V(G):g_{v_{i(u)}\to u}(x)\geq 1\}$ is dominating.
Note that $i(u)$ is well-defined, because $g$ is a weak Roman dominating function.
For $i\in [k]$, let $Q_i=\{ u\in V(G)\setminus (V_1\cup P_1\cup\cdots\cup P_k):i(u)=i\}$.
Since for every $u\in Q_i$, 
the set $\{ x\in V(G):g_{v_i\to u}(x)\geq 1\}$ is dominating,
we obtain that every vertex in $Q_i$ is adjacent to every vertex in $\{ v_i\}\cup P_i$,
which completes the proof. $\Box$

\begin{corollary}\label{corollary2}
Let $G$ be a connected $\{ K_4,K_4-e\}$-free graph.

$\gamma_{r2}(G)=2\gamma_r(G)$ if and only if 
\begin{itemize}
\item either $G$ is $K_2$,
\item or $G$ arises by adding a matching containing two edges between two disjoint triangles, 
\item or $G$ arises from the disjoint union of $k=\gamma_r(G)$ triangles
$v_1w_1u_1v_1$, $v_2w_2u_2v_2$, $\ldots$, $v_kw_ku_kv_k$
by adding edges between the vertices in $\{ v_1,\ldots,v_k\}$.
\end{itemize}
\end{corollary}
{\it Proof:} Since the sufficiency is straightforward, we only prove the necessity.
Therefore, let $G$ be a connected $\{ K_4,K_4-e\}$-free graph with $\gamma_{r2}(G)=2\gamma_r(G)$.
Let $g:V(G)\to \{ 0,1,2\}$ be a minimum weak Roman dominating function of $G$, 
and let $V_1,P_1,\ldots,P_k,Q_1,\ldots,Q_k$ be as in Theorem \ref{theorem1},
that is, $k=\gamma_r(G)$.
Since $G$ is $\{ K_4,K_4-e\}$-free, 
we have $|Q_i|\leq 1$ and $|P_i|+|Q_i|\leq 2$ for every $i\in [k]$.
This implies that $G$ has a spanning subgraph $H$ that is the union of $\ell$ triangles 
$v_1w_1u_1v_1$, 
$v_2w_2u_2v_2$,
$\ldots$,
$v_\ell w_\ell u_\ell v_\ell$
for some $\ell\leq k$, and 
$k-\ell$ complete graphs of order two
$v_{\ell+1}u_{\ell+1}$,
$v_{\ell+2}u_{\ell+2}$,
$\ldots$,
$v_ku_k.$

If a vertex $u'$ in some component $v_iu_i$ of $H$ with $\ell+1\leq i\leq k$ 
has a neighbor $v'$ in some other component $K$ of $H$, then 
$f:V(G)\to 2^{[2]}$ with
$$f(x)=
\left\{
\begin{array}{rl}
\{ 1,2\}, & x\in \{ v_1,\ldots,v_k\}\setminus (\{ v_i,u_i\}\cup V(K)),\\
\{ 1,2\}, & x=v',\\
\{ 1\}, & x\in \{ v_i,u_i\}\setminus \{ u'\},\mbox{ and }\\
\emptyset, & \mbox{ otherwise}
\end{array}
\right.
$$
is a $2$-rainbow dominating function of $G$ of weight less than $2\gamma_r(G)$,
which is a contradiction.
Since $G$ is connected, this implies that $G$ is either $K_2$ or $\ell=k$.
Hence, we may assume that $\ell=k$,
that is, $H$ is the union of $k$ triangles.

If there are two edges $v'u'$ and $v''u''$ such that $u',u''\in V(K)$ with $u'\not=u''$,
$v'\in V(K')$, and $v''\in V(K'')$ for three distinct components $K$, $K'$, and $K''$ of $H$,
then 
$f:V(G)\to 2^{[2]}$ with
$$f(x)=
\left\{
\begin{array}{rl}
\{ 1,2\}, & x\in \{ v_1,\ldots,v_k\}\setminus (V(K)\cup V(K')\cup V(K'')),\\
\{ 1,2\}, & x\in \{ v',v''\},\\
\{ 1\}, & x\in V(K)\setminus \{ u',u''\},\mbox{ and }\\
\emptyset, & \mbox{ otherwise}
\end{array}
\right.
$$
is a $2$-rainbow dominating function of $G$ of weight less than $2\gamma_r(G)$,
which is a contradiction.
Hence, such a pair of edges does not exist. 
In view of the desired statement, we may now assume that 
there is some component $K$ of $H$ such that two vertices in $K$ have neighbors in other components of $H$.
By the previous observation and since $G$ is connected, 
this implies that $k=2$, 
and that $G$ arises by adding a matching containing two or three edges between two disjoint triangles.
If $G$ arises by adding a matching containing three edges between two disjoint triangles,
then $\gamma_{r2}(G)=3<2\gamma_r(G)$, which is a contradiction.
Hence, $G$ arises by adding a matching containing two edges between two disjoint triangles,
which completes the proof. $\Box$

\bigskip

\noindent The last result immediately implies the following.

\begin{corollary}\label{corollary1}
Let $G$ be a triangle-free graph.

$\gamma_{r2}(G)=2\gamma_r(G)$ if and only if 
$G$ is the disjoint union of copies of $K_2$.
\end{corollary}

\begin{theorem}\label{theorem2}
It is NP-hard to decide $\gamma_{r2}(G)=2\gamma_r(G)$ for a given $K_5$-free graph $G$.
\end{theorem}
{\it Proof:} We describe a reduction from {\sc 3Sat}.
Therefore, let $F$ be a {\sc 3Sat} instance
with $m$ clauses $C_1,\ldots,C_m$
over $n$ boolean variables $x_1,\ldots,x_n$.
Clearly, we may assume that $m\geq 2$.
We will construct a $K_5$-free graph $G$ whose order is polynomially bounded in terms of $n$ and $m$
such that $F$ is satisfiable if and only if $\gamma_{r2}(G)=2\gamma_r(G)$.
For every variable $x_i$, create a copy $G(x_i)$ of $K_4$ 
and denote two distinct vertices of $G(x_i)$ by $x_i$ and $\bar{x}_i$.
For every clause $C_j$, create a vertex $c_j$.
For every literal $x\in \{ x_1,\ldots,x_n\}\cup \{ \bar{x}_1,\ldots,\bar{x}_k\}$ 
and every clause $C_j$ such that $x$ appears in $C_j$,
add the edge $xc_j$.
Finally, add 
two further vertices $a$ and $b$, 
the edge $ab$, 
and all possible edges between $\{ a,b\}$ and $\{ c_1,\ldots,c_m\}$.
This completes the construction of $G$.
Clearly, $G$ is $K_5$-free and has order $4n+m+2$.

Let $f$ be a $2$-rainbow dominating function of $G$.
Clearly, $\sum_{u\in \{ a,b\}\cup \{ c_1,\ldots,c_m\}}|f(u)|\geq 2$, and 
$\sum_{u\in V(G_i)}|f(u)|\geq 2$ for every $i\in [n]$, 
which implies $\gamma_{r2}(G)\geq 2n+2$.
Since 
$$x\mapsto
\left\{
\begin{array}{rl}
\{ 1,2\}, & x\in \{ a,x_1,\ldots,x_n\}\mbox{ and }\\
\emptyset, & \mbox{ otherwise}
\end{array}
\right.
$$
defines a $2$-rainbow dominating function of weight $2n+2$,
we obtain $\gamma_{r2}(G)=2n+2$.
By Theorem \ref{theorem1}, we have $\gamma_r(G)\geq n+1$.
In remains to show that $F$ is satisfiable if and only if $\gamma_r(G)=n+1$.

Let $\gamma_r(G)=n+1$.
Let $g$ be a minimum weak Roman dominating function of $G$.
By Theorem \ref{theorem1}, there is no vertex $x$ of $G$ with $g(x)=2$.
Let $V_1$ be the set of vertices $x$ of $G$ with $g(x)=1$.
Since, $\sum_{u\in \{ a,b\}\cup \{ c_1,\ldots,c_m\}}g(u)\geq 1$, and 
$\sum_{u\in V(G_i)}g(u)\geq 1$ for every $i\in [n]$, 
we obtain that $\{ a,b\}\cup \{ c_1,\ldots,c_m\}$ contains exactly one vertex, say $y_0$, from $V_1$,
and that $V(G_i)$ contains exactly one vertex, say $y_i$, from $V_1$ for every $i\in [n]$.
Since $m\geq 2$, we may assume, by symmetry, that $g(c_1)=0$.
If no neighbor $v$ of $c_1$ with $g(v)\geq 1$ belongs to $\{ a,b\}\cup \{ c_1,\ldots,c_m\}$, 
then $g$ is not a weak Roman dominating function.
Hence $y_0\in \{ a,b\}$, and $y_0$ is the only neighbor of $c_1$ with positive $g$-value
such that the set $\{ x\in V(G):g_{y_0\to c_1}(x)\geq 1\}$ is dominating,
which implies that for every $\ell\in [m]\setminus \{ 1\}$,
the vertex $c_\ell$ is adjacent to a vertex in $\{ y_1,\ldots,y_k\}$.
Since $y_0\in \{ a,b\}$ and $m\geq 2$, this actually implies, by symmetry, that for every $\ell\in [m]$,
the vertex $c_\ell$ is adjacent to a vertex in $\{ y_1,\ldots,y_k\}$,
that is, the intersection of $\{ y_1,\ldots,y_k\}$
with $\{ x_1,\ldots,x_n\}\cup \{ \bar{x}_1,\ldots,\bar{x}_k\}$
indicates a satisfying truth assignment for $F$.

Conversely, if $F$ has a satisfying truth assignment, then
$$x\mapsto
\left\{
\begin{array}{rl}
1, & x=a,\\
1, & x\in \{ x_1,\ldots,x_n\}\cup \{ \bar{x}_1,\ldots,\bar{x}_k\}\mbox{ and $x$ is true, and}\\
0, & \mbox{ otherwise}
\end{array}
\right.
$$
defines a weak Roman dominating function of $G$ of weight $n+1$,
which implies $\gamma_r(G)=n+1$, and completes the proof. $\Box$

\bigskip

\noindent For a positive integer $k$, let 
$${\cal G}_k=\{ G:\forall H\subseteq_{\rm ind}G:\gamma_r(H)\geq k\Rightarrow \gamma_{r2}(H)=2\gamma_r(H)\},$$
where $H\subseteq_{\rm ind}G$ means that $H$ is an induced subgraph of $G$.
Since $\gamma_{r2}(K_1)=1=\gamma_r(K_1)$,
the set ${\cal G}_1$ contains no graph of positive order.
Since $\gamma_{r2}(\bar{K}_2)=2=\gamma_r(\bar{K}_2)$,
where $\bar{H}$ denotes the complement of some graph $H$,
the set ${\cal G}_2$ consists exactly of all complete graphs.
The smallest value for $k$ that leads to an interesting class of graphs is $3$.

\begin{theorem}\label{theorem3} 
${\cal G}_3={\rm Free}(\{ \bar{K}_3,C_5\})$.
\end{theorem} 
{\it Proof:} 
Since $\gamma_{r2}(\bar{K}_3)=3=\gamma_r(\bar{K}_3)$ and 
$\gamma_{r2}(C_5)=3=\gamma_r(C_5)$, it follows easily that 
$\bar{K}_3$ and $C_5$ are minimal forbidden induced subgraphs for ${\cal G}_3$.
Now, let $G$ be a minimal forbidden induced subgraphs for ${\cal G}_3$,
which implies that $\gamma_r(G)\geq 3$ and $\gamma_{r2}(H)\not=2\gamma_r(H)$.
It remains to show that $G$ is either $\bar{K}_3$ or $C_5$.
For a contradiction, we assume that $G$ is neither $\bar{K}_3$ nor $C_5$.
Since $G$ is a minimal forbidden induced subgraph, this implies that $G$ is $\{ \bar{K}_3,C_5\}$-free.
Since $\gamma_r(G)\geq 3$, the graph $G$ is not complete.
Let $u$ and $v$ be two non-adjacent vertices of $G$.
Since $G$ is $\bar{K}_3$-free, 
we have $V(G)= \{ u,v\}\cup N_u\cup N_v\cup N_{u,v}$, 
where 
$N_u=N_G(u)\setminus N_G(v)$,
$N_v=N_G(v)\setminus N_G(u)$, and
$N_{u,v}=N_G(u)\cap N_G(v)$.
Since $G$ is $\bar{K}_3$-free, 
the sets $N_u$ and $N_v$ are complete.
If for every vertex $w$ in $N_{u,v}$, we have $N_u\subseteq N_G(w)$ or $N_v\subseteq N_G(w)$,
then 
$$x\mapsto
\left\{
\begin{array}{rl}
1, & x\in \{ u,v\},\mbox{ and}\\
0, & \mbox{ otherwise}
\end{array}
\right.
$$
defines a weak Roman dominating function of $G$ of weight $2$,
which implies the contradiction $\gamma_r(G)<3$.
Hence, there are vertices $w_u\in N_u$, $w_v\in N_v$, and $w_{u,v}\in N_{u,v}$
such that $w_{u,v}$ is adjacent to neither $w_u$ nor $w_v$.
Since $G$ is $\bar{K}_3$-free, this implies that $w_u$ is adjacent to $w_v$, 
and $uw_uw_vvw_{u,v}u$ is an induced $C_5$ in $G$,
which is a contradiction and completes the proof. $\Box$

\bigskip

\noindent Our results motivate several questions.
Do the graphs $G$ with $\gamma_{r2}(G)=2\gamma_r(G)$
that are either $K_4$-free or $(K_4-e)$-free have a simple structure?
Can they at least be recognized efficiently?
Can Theorem \ref{theorem2} be strengthened by restricting the input graphs even further?
What are the minimal forbidden induced subgraphs for the classes ${\cal G}_k$ where $k\geq 4$?

\bigskip

\noindent {\bf Acknowledgment}  
J.D. Alvarado and S. Dantas were partially supported by FAPERJ, CNPq, and CAPES.


\begin{thebibliography}{}
\bibitem{bhr} B. Bre\v{s}ar, M.A. Henning, and D.F. Rall, Rainbow domination in graphs, Taiwanese J. Math. 12 (2008) 213-225.
\bibitem{chh} M. Chellali, T.W. Haynes, and S.T. Hedetniemi, Bounds on weak roman and 2-rainbow domination numbers, Discrete Appl. Math. 178 (2014) 27-32.
\bibitem{cr} M. Chellali and N.J. Rad, On 2-rainbow domination and Roman domination in graphs, Australas. J. Combin. 56 (2013) 85-93.
\bibitem{ff} S. Fujita and M. Furuya, Difference between 2-rainbow domination and Roman domination in graphs, Discrete Appl. Math. 161 (2013) 806-812.
\bibitem{hh} M.A. Henning and S.T. Hedetniemi, Defending the Roman Empire - new strategy, Discrete Math. 266 (2003) 239-251.
\bibitem{st} I. Stewart, Defend the Roman empire!, Sci. Am. 281 (1999) 136-139.
\bibitem{wx} Y. Wu and H. Xing, Note on 2-rainbow domination and Roman domination in graphs, Appl. Math. Lett. 23 (2010) 706-709.
\end{thebibliography}
\end{document}